# Geometric Interpretation of $\det(\alpha^{(3)}, \alpha^{(4)}, \alpha^{(5)})$ in $E_1^3$

Seher KAYA, İsmail GÖK and Yusuf YAYLI

ABSTRACT. In this paper, we investigate the tangent indicatrix of the curve $\alpha$ with constant curvature. Tangent indicatrix of the curve $\alpha$ is characterized with $\det(\alpha^{(3)}, \alpha^{(4)}, \alpha^{(5)}) = 0$ in Minkowski 3-space $E_1^3$. Moreover, we study null slant helices using the determinant approach and give the following characterization:
  A curve $\alpha$ is a null slant helix in $E_1^3$ if and only if $\det(\alpha^{(3)}, \alpha^{(4)}, \alpha^{(5)}) = 0$.
Then similar results are obtained for non-null curves with the condition $\kappa = 1$.

# Introduction

Let $\alpha = \alpha(s)$ in Minkowski $3-$ space with the Frenet apparatus $\{T, N, B, \kappa, \tau\}$, where $\kappa$ and $\tau$ represent curvature and torsion of the curve $\alpha$. A general helix in Euclidean $3-$ space is a curve whose tangent vector makes a constant angle with a fixed direction (the axis of the general helix). By using its curvature functions we have a necessary and sufficient condition that the ratio of curvature to torsion be constant. A slant helix defined by the property that its normal vector makes a constant angle with a fixed straight line (the axis of the slant helix) in Euclidean $3-$ space was introduced by Izumiya and Takeuchi. Also, they proved that $\alpha$ is a slant helix if and only if $\dfrac{\kappa^2}{(\kappa^2 + \tau^2)^{\frac{3}{2}}} \left(\dfrac{\tau}{\kappa}\right)'$ is a constant function.

In $[10]$, Takenaka considered curves characterized by $\det(\alpha^{(3)}, \alpha^{(4)}, \alpha^{(5)}) = 0$ in Euclidean $3-$ space and gave the following theorem:
(1) The condition $\det(\alpha^{(0)}, \alpha^{(1)}, \alpha^{(2)}) = 0$ characterizes a great circle, where $\alpha$ is a spherical curve.
(2) The condition $\det(\alpha^{(1)}, \alpha^{(2)}, \alpha^{(3)}) = 0$ characterizes a plane curve.
(3) The condition $\det(\alpha^{(2)}, \alpha^{(3)}, \alpha^{(4)}) = 0$ characterizes a curve of constant slope, where $\alpha^{(k)}$ denote the $k$-th differential by arc-lenght parameter of the curve $\alpha$.
After this study, Yaylı and Saracoğlu $[11]$ investigated the curves $\alpha$ with the condition $\det(\alpha^{(3)}, \alpha^{(4)}, \alpha^{(5)}) = 0$ and proved that the curves are Salkowski curves.
Ali and Lopez $[1]$ gave characterizations of a slant helix $\alpha = \alpha(s)$ in terms of its curvature and torsion in Minkowski $3-$ space. The authors investigated the tangent and binormal indicatrices of the slant helix $\alpha$ and proved that they are helices in $E_1^3$. Also, Gök et all. gave some differential equations which characterize timelike slant helices by using their Frenet apparatus of spherical indicatrices $[5]$.
Salkowski curves defined by E. Salkowski which have a constant curvatures and non-constant torsion are special family of the slant helices [9]. In Minkowski $3-$ space, spacelike and timelike Salkowski curves were investigated [2-4]. The following figures are examples of Salkowski curves in the paper [3].

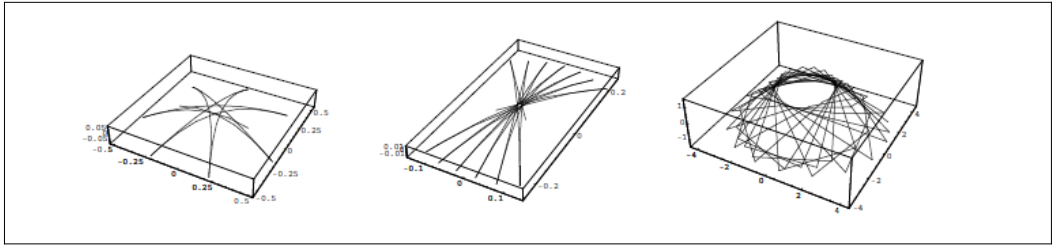

Fig 1. Some Spacelike Salkowski Curves

In this paper, we characterize the tangent indicatrix $(T)$ of the curve $\alpha$ with the condition $\kappa = 1$ and $\det(\alpha^{(3)}, \alpha^{(4)}, \alpha^{(5)}) = 0$ in $E_1^3$. Then the relationship between the determinant function and non-null slant curves with the condition $\kappa = 1$ is given. Furthermore, we show that a curve is a lightlike slant helix if and only if $\det(\alpha^{(3)}, \alpha^{(4)}, \alpha^{(5)}) = 0$. According to our idea this characterization is useful for the null slant helix. Moreover, we verify the torsion of the Salkowski curves considering the condition $\det(\alpha^{(3)}, \alpha^{(4)}, \alpha^{(5)}) = 0$.

## Preliminaries

In this section, we review some basic concepts on classical differential geometry of space curves in Minkowski $3-$ space.

For any two vectors $x = (x_1, x_2, x_3)$ and $y = (y_1, y_2, y_3) \in E_1^3$, Lorentzian inner product of $x$ and $y$ as the following

$$g(x, y) = -x_1 y_1 + x_2 y_2 + x_3 y_3$$

Let $u = (u_1, u_2, u_3)$ be a vector in Minkowski $3-$ space $E_1^3$. The vector $u$ is called spacelike if $g(u,u) > 0$ or $u = 0$, timelike if $g(u,u) < 0$ and lightlike if $g(u,u) = 0$. The norm of the vector $u$ is given by $\|u\| = \sqrt{|g(u,u)|}$. An arbitrary curve $\alpha : I \subset \mathbb{R} \to E_1^3$ is called spacelike (timelike or lightlike), if its velocity vector $\alpha'(s)$ is spacelike (timelike or lightlike). A non-null curve $\alpha$ is arc-lenght parametrized if $g(\alpha'(s), \alpha'(s)) = 1$ and a null curve $\alpha$ is parametrized by pseudo-arc $s$ if $g(\alpha''(s), \alpha''(s)) = 1$ [8].

Let $\alpha : I \subset \mathbb{R} \to E_1^3$ be a non-null unit speed space curve in $E_1^3$ and $\{T, N, B\}$ denotes the Frenet frame of the curve $\alpha$. Then the Frenet equations are given

(1)
$$T' = \kappa N$$
$$N' = -\varepsilon_T \varepsilon_N \kappa T + \tau B$$
$$B' = -\varepsilon_N \varepsilon_B \tau N$$

where $\varepsilon_T$, $\varepsilon_N$, $\varepsilon_B$ are the causal characters of tangent vector, normal vector and binormal vector of the curve $\alpha$ respectively [7].

For a null space curve $\alpha(s)$ with $\kappa$ and $\tau$ curvature functions in the space $E_1^3$, the following Frenet formulae are given in [8]

(2)
$$T' = \kappa N$$
$$N' = \tau T - \kappa B$$
$$B' = -\tau N$$

$g(T,T) = g(B,B) = 0, \quad g(N,N) = 1, \quad g(T,N) = g(N,B) = 0, \quad g(T,B) = 1.$

In this case, $\kappa$ can be two value:
i) $\kappa = 0$, $\alpha$ is a straight null line
ii) $\kappa = 1$, all other cases.

# Tangents Indicatrix of Non-Null Curves with Constant Curvature

In this section, we give some characterizations of the tangent indicatrix of the non-null curve by using $\det(\alpha^{(3)}, \alpha^{(4)}, \alpha^{(5)})$.

**Proposition 1.** *Let* $\alpha = \alpha(s)$ *be a unit speed non-null curve with the condition* $\kappa = 1$ *in Minkowski* $3-$ *space. The tangent indicatrix of the space curve* $\alpha$ *is spherical helix if and only if* $\det(\alpha^{(3)}, \alpha^{(4)}, \alpha^{(5)}) = 0$.

**Proof.** Let $\alpha$ be a non-null curve with the condition $\kappa = 1$ in $E_1^3$ and $(T)$ be the tangent indicatrix of the curve $\alpha$. Then

(3)
$$\alpha'(s) = (T)$$

if we differentiate the Eq. (3), we have
$$T' = \kappa N = N$$

Since $\|T'\| = 1$, $(T)$ is a unit-speed curve. The curve $\alpha$ and its tangent indicatrix $(T)$ have the same parameter. If $(T)$ is helix, then
$$\det(T^{(2)}, T^{(3)}, T^{(4)}) = 0$$

By using the equality $T = \alpha'$, we obtain
$$\det(\alpha^{(3)}, \alpha^{(4)}, \alpha^{(5)}) = 0.$$

Conversely, we assume that $\det(\alpha^{(3)}, \alpha^{(4)}, \alpha^{(5)})$ is zero function. We can easily see that $(T)$ is a spherical helix. This completes the proof.

**Remark 1.** *If the curve* $\alpha$ *is a spacelike curve with lightlike normal vector in* $E_1^3$, *then the tangent indicatrix of the curve* $\alpha$ *is a spherical helix.*

**Proposition 2.** *Let* $\alpha = \alpha(s)$ *be a unit speed non-null curve with the condition* $\kappa = 1$ *in* $E_1^3$. *The tangent indicatrix* $(T)$ *of the curve* $\alpha$ *is spherical helix if and only if*

(4)
$$\tau''(1 + \tau^2 \varepsilon_T \varepsilon_B) - 3\tau(\tau')^2 \varepsilon_T \varepsilon_B = 0$$

*where* $\varepsilon_T$ *and* $\varepsilon_B$ *are the causal characters of the tangent and binormal vectors of the curve* $\alpha$, *respectively.*

**Proof.** Let the tangent indicatrix of the curve $\alpha$ be a spherical helix. By using Prop. 1, we have $\det(\alpha^{(3)}, \alpha^{(4)}, \alpha^{(5)}) = 0$. If we calculate the derivative of the curve $\alpha$ up to the fifth order, we get following equalities given by

$$\alpha' = T$$
$$\alpha'' = \kappa N = N, \quad \kappa = 1$$
$$\alpha^{(3)} = -\varepsilon_T \varepsilon_N T + \tau B$$
$$\alpha^{(4)} = (-\varepsilon_T \varepsilon_N - \tau^2 \varepsilon_N \varepsilon_B) N + \tau' B$$
$$\alpha^{(5)} = (1 + \tau^2 \varepsilon_T \varepsilon_B) T + (-3\tau\tau' \varepsilon_N \varepsilon_B) N + (-\tau \varepsilon_T \varepsilon_N - \tau^3 \varepsilon_N \varepsilon_B + \tau'') B$$

Since $\det(\alpha^{(3)}, \alpha^{(4)}, \alpha^{(5)}) = 0$, we obtain the Eq. (4).

Conversely, assume that the Eq. (4) is provided. It is obvious that the tangent indicatrix of the curve $\alpha$ is a spherical helix. Which completes the proof.

**Corollary 1.** *Let* $\alpha : I \subset \mathbb{R} \to E_1^3$ *be a unit speed curve with the condition* $\kappa = 1$ *in* $E_1^3$.
*Case I: If* $\alpha$ *is a spacelike curve with spacelike normal vector or a timelike curve, the tangent indicatrix of the space curve* $\alpha$ *is a spherical helix if and only if the equality*

(5) $$\tau''(1 - \tau^2) + 3\tau(\tau')^2 = 0$$

*holds.*

*Case II: If* $\alpha$ *is a spacelike curve with timelike normal vector field, the tangent indicatrix of the space curve* $\alpha$ *is a spherical helix if and only if the equality*

(6) $$\tau''(1 + \tau^2) - 3\tau(\tau')^2 = 0$$

*holds.*

Then, using the solution of the Eq. (5) and Eq. (6) we have following corollary.

**Corollary 2.** *Let* $\alpha : I \subset \mathbb{R} \to E_1^3$ *be a unit speed curve with the condition* $\kappa = 1$ *in* $E_1^3$.
*Case I: If* $\alpha$ *is spacelike curve with spacelike normal vector field or a timelike curve. The tangent indicatrix of the space curve* $\alpha$ *is a spherical helix if and only if the torsion of the curve* $\alpha$ *is*

(7) $$\tau = \pm \frac{bs + c}{\left[\pm 1 + (bs + c)^2\right]^{\frac{1}{2}}}$$

*where*

$$\frac{1}{b}(1 - c) \leq s \leq -\frac{1}{b}(1 + c), \quad b \neq 0 \text{ and } b, c \in \mathbb{R}$$

*Case II: If* $\alpha$ *is a spacelike curve with timelike normal vector field. The tangent indicatrix of the space curve* $\alpha$ *is a spherical helix if and only if the torsion of the curve* $\alpha$ *is*

(8) $$\tau = \pm \frac{bs + c}{\left[1 - (bs + c)^2\right]^{\frac{1}{2}}}$$

*where*

$$-\frac{1}{b}(1 + c) \leq s \leq \frac{1}{b}(1 - c), \quad b \neq 0 \text{ and } b, c \in \mathbb{R}$$

**Proof.** Let $\alpha : I \subset \mathbb{R} \to E_1^3$ be a unit speed spacelike curve which has a spacelike normal vector field or timelike curve with the condition $\kappa = 1$ in $E_1^3$ and the tangent indicatrix of the space curve $\alpha$ be spherical helix. By using Corollary 1 we have the Eq. (5). Then, we can easily see that the solutions of the differential equation given Eq. (5) are

$$\tau = \pm \frac{bs+c}{\left[\pm 1 + (bs+c)^2\right]^{\frac{1}{2}}}.$$

Conversely, it is easy to obtain that the Eq. (7) holds in the Eq. (5). From Corollary 1, proof is completed.

Similarly, Case II can be proved easily.

**Theorem 1.** *Let $\alpha$ be a unit speed spacelike curve in $E_1^3$.*
*(a) If the normal vector of $\alpha$ is spacelike, then $\alpha$ is a slant helix if and only if either one the next two functions*

$$\frac{\kappa^2}{(\tau^2 - \kappa^2)^{\frac{3}{2}}} \left(\frac{\tau}{\kappa}\right)' \quad \text{or} \quad \frac{\kappa^2}{(\kappa^2 - \tau^2)^{\frac{3}{2}}} \left(\frac{\tau}{\kappa}\right)'$$

*is constant everywhere $\tau^2 - \kappa^2$ does not vanish.*
*(b) If the normal vector of $\alpha$ is timelike, then $\alpha$ is a slant helix if and only if the function*

$$\frac{\kappa^2}{(\tau^2 + \kappa^2)^{\frac{3}{2}}} \left(\frac{\tau}{\kappa}\right)'$$

*is constant.*
*(c) Any spacelike curve with lightlike normal vector is a slant curve* [1].

**Theorem 2.** *Let $\alpha$ be a unit speed timelike curve in $E_1^3$. Then $\alpha$ is a slant helix if and only if either one the next two functions*

$$\frac{\kappa^2}{(\tau^2 - \kappa^2)^{\frac{3}{2}}} \left(\frac{\tau}{\kappa}\right)' \quad \text{or} \quad \frac{\kappa^2}{(\kappa^2 - \tau^2)^{\frac{3}{2}}} \left(\frac{\tau}{\kappa}\right)'$$

*is constant everywhere $\tau^2 - \kappa^2$ does not vanish* [1].

**Proposition 3.** *Let $\alpha$ be a unit speed non-null curve with the condition $\kappa = 1$. The curve $\alpha$ is a timelike or spacelike (with spacelike or timelike normal vector field) slant helix if and only if $\det(\alpha^{(3)}, \alpha^{(4)}, \alpha^{(5)}) = 0$.*

**Proof.** Let $\alpha$ be a unit speed timelike or spacelike slant helix with spacelike normal vector field, in that case the tangent indicatrix $(T)$ of the curve $\alpha$ is a spherical helix [1]. By using the Prop.1, we have $\det(\alpha^{(3)}, \alpha^{(4)}, \alpha^{(5)}) = 0$.

On the contrary, if $\det(\alpha^{(3)}, \alpha^{(4)}, \alpha^{(5)}) = 0$ we obtain the Eq. (5). Also if we calculate the derivative of the equation

$$\frac{\kappa^2}{(\tau^2 - \kappa^2)^{\frac{3}{2}}} \left(\frac{\tau}{\kappa}\right)' = \text{constant}$$

we get

$$\tau''(1 - \tau^2) + 3\tau(\tau')^2 = 0$$

where $\kappa = 1$.

So, we can easily see that $\det(\alpha^{(3)}, \alpha^{(4)}, \alpha^{(5)}) = 0$. Then, we have $\left(\frac{\kappa^2}{(\tau^2 - \kappa^2)^{\frac{3}{2}}} \left(\frac{\tau}{\kappa}\right)'\right)' = 0.$

That is $\frac{\kappa^2}{(\tau^2-\kappa^2)^{\frac{3}{2}}}(\frac{\tau}{\kappa})' =$ constant. From Theorem 1, we say that $\alpha$ is a spacelike (with spacelike normal vector field) slant helix. Also the same case is valid for timelike curve. Similarly we can show that $\alpha$ is a spacelike slant helix with timelike normal vector field if and only if $\det(\alpha^{(3)}, \alpha^{(4)}, \alpha^{(5)}) = 0$.

**Corollary 3.** *Let $\alpha$ be a unit speed spacelike curve with the condition $\kappa = 1$ in $E_1^3$. Then we have the following cases:*

*a) The curve $\alpha$ is a spacelike slant with spacelike normal vector field or $\alpha$ is a timelike slant helix if and only if the torsion of the curve $\alpha$ is $\tau = \pm \dfrac{bs+c}{\left[\pm 1 + (bs+c)^2\right]^{\frac{1}{2}}}$.*

*b) The curve $\alpha$ is a spacelike slant helix with timelike normal vector field if and only if the torsion of the curve $\alpha$ is $\tau = \pm \dfrac{bs+c}{\left[1 - (bs+c)^2\right]^{\frac{1}{2}}}$*

**Corollary 4.** *Let $\alpha = \alpha(s)$ be a unit speed curve with the condition $\kappa = 1$ in Minkowski space $E_1^3$. The following three conditions are equivalent:*

1) The space curve $\alpha$ is a spacelike (has a spacelike normal vector field) or a timelike Salkowski slant curve
2) $\det(\alpha^{(3)}, \alpha^{(4)}, \alpha^{(5)}) = 0$
3) $\tau = \pm \dfrac{bs+c}{\left[\pm 1 + (bs+c)^2\right]^{\frac{1}{2}}}$

**Proof.** By using Prop.3, we say that if the space curve $\alpha$ is a spacelike (with spacelike normal vector field) or a timelike slant helix, then $\det(\alpha^{(3)}, \alpha^{(4)}, \alpha^{(5)}) = 0$. From the $\det(\alpha^{(3)}, \alpha^{(4)}, \alpha^{(5)}) = 0$ we obtain the Eq. (5). The solutions of the equation are

$$\tau = \pm \frac{bs+c}{\left[\pm 1 + (bs+c)^2\right]^{\frac{1}{2}}}.$$

**Corollary 5.** *Let $\alpha = \alpha(s)$ be a unit speed curve with the condition $\kappa = 1$ in Minkowski space $E_1^3$. The following three conditions are equivalent:*

1) The space curve $\alpha$ is spacelike Salkowski slant curve with timelike normal vector field.
2) $\det(\alpha^{(3)}, \alpha^{(4)}, \alpha^{(5)}) = 0$.
3) $\tau = \pm \dfrac{bs+c}{\left[1 - (bs+c)^2\right]^{\frac{1}{2}}}$.

**Remark 2.** *In fact, Ali gave the torsions of the spacelike and timelike Salkowski curves with the condition $\kappa = 1$ in $[2, 6, 7]$. If we take $c = 0$ and $b = \dfrac{1}{\tanh\phi}$ in Corollary 4 and Corollary 5, then we verify these statement with the help of determinant approach.*

*Case I: If the curve $\alpha$ is a spacelike (with spacelike normal vector field) Salkowski curve we get*

$$\tau = \pm \frac{s}{\left[\tanh^2 \phi + s^2\right]^{\frac{1}{2}}}.$$

**Case II:** *If the curve* $\alpha$ *is a spacelike (with timelike normal vector field) Salkowski curve we get*

$$\tau = \pm \frac{s}{\left[\tanh^2 \phi - s^2\right]^{\frac{1}{2}}}.$$

**Case III:** *If the curve* $\alpha$ *is a timelike Salkowski curve we get* $\tau = \pm \dfrac{s}{\left[s^2 - \tanh^2 \phi\right]^{\frac{1}{2}}}.$

# Characterization of Null Slant Helix with $\det(\alpha^{(3)},\alpha^{(4)},\alpha^{(5)})=0$

In this section, we will characterize null slant helices in Minkowski space $E_1^3$ using the determinant approach. If the curve $\alpha$ is not a straight null line it has the curvature $\kappa=1$. So, we can examine the curves with the help of our determinant approach.

**Theorem 3.** *Let* $\alpha=\alpha(s)$ *be a unit speed lightlike curve in Minkowski space* $E_1^3$. *Then* $\alpha$ *is a slant helix if and only if the torsion of the curve is*

$$\tau(s)=\frac{a}{(bs+c)^2}$$

*where* $a,b$ *and* $c$ *are constants, with* $bs+c\neq 0$ [1].

**Proposition 4.** *Let* $\alpha=\alpha(s)$ *be a unit speed lightlike curve in Minkowski space* $E_1^3$. *The curve* $\alpha$ *is a null ( lightlike ) slant helix if and only if*
$$\det(\alpha^{(3)},\alpha^{(4)},\alpha^{(5)})=0.$$

**Proof.** Let the curve $\alpha=\alpha(s)$ be a lightlike slant helix. From Theorem 3, the torsion of the curve $\alpha$ is $\tau(s)=\frac{a}{(bs+c)^2}$. $\tau$ is the solution of the differentiable equation

$$2\tau\tau''-3(\tau')^2=0.$$

Since $\det(\alpha^{(3)},\alpha^{(4)},\alpha^{(5)})=2\tau\tau''-3(\tau')^2$, we can easily say that $\det(\alpha^{(3)},\alpha^{(4)},\alpha^{(5)})=0$.
Let $\alpha$ be a lightlike curve in $E_1^3$ with the condition $\det(\alpha^{(3)},\alpha^{(4)},\alpha^{(5)})=0$. From the Frenet equations for lightlike curve given Eq.(1), we can calculate derivatives of the curve $\alpha$ as the following

$$\alpha^{(3)}=\tau T-B$$
$$\alpha^{(4)}=\tau' T+2\tau N$$
$$\alpha^{(5)}=(\tau''+2\tau^2)T+3\tau' N-2\tau B.$$

and then the equality $\det(\alpha^{(3)},\alpha^{(4)},\alpha^{(5)})=0$ gives us
$$2\tau\tau''-3(\tau')^2=0.$$

The solution of the last equation is $\tau(s)=\frac{a}{(bs+c)^2}$. So, using the Theorem 2, we can easily obtain that $\alpha$ is a lightlike slant helix. Which completes the proof.

**Example1.** *Let* $\alpha(s) = \frac{1}{6}(\frac{s^5}{5} - \frac{1}{s}, s^2, \frac{s^5}{5} + \frac{1}{s})$ *be a null curve in* $E_1^3$. *The curve* $\alpha$ *is null W- slant helix with Cartan curvature functions* $k_1 = 1$ *and* $k_2 = -\frac{4}{s^2}$ [6].

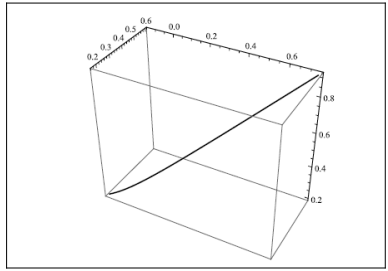

Fig 2. Null W-slant helix

*Since* $\det(\alpha^{(3)}, \alpha^{(4)}, \alpha^{(5)}) = 0$, *by using the Prop.* 4 *we can easily see that* $\alpha$ *is a null slant helix.*

# Conclusion

Characterizations of the curves are very important in terms of differential geometry. In this paper, we have characterized non-null slant helices with the condition $\kappa = 1$ and Salkowski curves with a new approach in Minkowski 3-space $E_1^3$. Also, a useful way is given for characterization of null slant helices.

Department of Mathematics, Faculty of Science, University of Ankara Tandogan, Ankara, TURKEY

E-mail address: seherkaya@ankara.edu.tr

Department of Mathematics, Faculty of Science, University of Ankara Tandogan, Ankara, TURKEY

E-mail address: igök@ankara.edu.tr

Department of Mathematics, Faculty of Science, University of Ankara Tandogan, Ankara, TURKEY

E-mail address: yyayli@ankara.edu.tr